# SMALL NOISE ASYMPTOTIC OF THE TIMING JITTER IN SOLITON TRANSMISSION

By Arnaud Debussche and Eric Gautier

*ENS Cachan Bretagne*

We consider the problem of the error in soliton transmission in long-haul optical fibers caused by the spontaneous emission of noise inherent to amplification. We study two types of noises driving the stochastic focusing cubic one dimensional nonlinear Schrödinger equation which appears in physics in that context. We focus on the fluctuations of the mass and arrival time or timing jitter. We give the small noise asymptotic of the tails of these two quantities for the two types of noises. We are then able to prove several results from physics among which the Gordon–Haus effect which states that the fluctuation of the arrival time is a much more limiting factor than the fluctuation of the mass. The physical results had been obtained with arguments difficult to fully justify mathematically.

**1. Introduction.** The nonlinear Schrödinger (NLS) equation occurs as a generic model in many areas of physics and describes the propagation of slowly varying envelopes of a wave packet in media with both weakly nonlinear and dispersive responses; see [35] for a detailed presentation. The one-dimensional equation with a cubic focusing nonlinearity, for example, has the form

$$i\frac{\partial u^{u_0}}{\partial t} = \Delta u^{u_0} + |u^{u_0}|^2 u^{u_0}, \tag{1.1}$$

where $u^{u_0}$ is a complex valued function depending on $t \geq 0$ and $x \in \mathbb{R}$ and the superscript $u_0$ means that $u^{u_0}(0,x) = u_0(x), x \in \mathbb{R}$. This equation is a very accurate model in the context of single-mode optical fibers over short distances. A derivation of the equation in that context is given, for example, in [24]. Resulting from a balance between the focusing nonlinearity and the dispersive linear part, localized stationary waves propagate. They are called solitons and have the form $\sqrt{2}A\,\mathrm{sech}(A(x-x_0))\exp(-iA^2t+i\theta_0)$, where $A>0$ is the amplitude, $x_0$ and $\theta_0$ are respectively the initial position and phase. By extension, we herein also call solitons the following nonstationary progressive solutions:

$$\sqrt{2}A\,\mathrm{sech}(A(x-x_0)+2AVt)\exp(-i(A^2-V^2)t+iV(x-x_0)+i\theta_0), \tag{1.2}$$

where $V$ is the group velocity or angular carrier frequency.

In soliton based amplitude-shifted-keyed (ASK) communication systems, solitons are used as information carriers to transmit the datum 0 or 1. A 1 corresponds to the emission of a soliton at time 0 with zero velocity $\Psi_A^0(x) = \sqrt{2}A\,\mathrm{sech}(Ax)$. It is produced by a laser beam. At coordinate $T$ (end of the line) a receiver records

$$(1/l)\int_{-l/2}^{l/2}|u^{u_0}(T,x)|^2\,dx, \qquad u_0=0 \quad \text{or} \quad u_0=\Psi_A^0.$$

In optics the usual $x$ variable of the NLS equation denotes some retarded time while $t$ is space. Thus, $[-l/2,l/2]$ is a window in time and $l$ may be chosen small since the solution $u^{u_0}$ of the NLS equation is localized and remains centered. When the above quantity is above a threshold, it is decided that a 1 has been emitted, otherwise it is decided that a 0 has been emitted.

Over long distances, damping induced low loss becomes significant and the signal has to be amplified. However, due to quantum considerations, amplification is intrinsically associated with small noise; see [15] for a physical justification of noise in optics. Simply stated, due to the Heisenberg principle, there is inherently uncertainty on the amplified signal. This uncertainty is accounted by noise in the system. This phenomenon is called spontaneous emission of noise. These intrinsic quantum features can have direct macroscopic consequences such as fluctuation of the arrival time, also called timing jitter or diffusion of the soliton. A practical consequence for engineering is error in soliton transmission.

We consider in this article two specific models from the physics literature on the topic which are stochastic PDEs (SPDEs). A first type of amplification, and most discussed, is the case of regularly spaced Erbium–Doped amplifiers placed along the line and such that the distance between amplifiers is small compared to the length of the line. The limit case where there is a continuum of amplifiers is called distributed amplification. In that case the noise (Gordon–Haus noise) acts as a random external force; see, for example, [14, 19, 31]. There, the following equation is used:

$$i\frac{\partial u^{\varepsilon,u_0}}{\partial t} = \Delta u^{\varepsilon,u_0} + |u^{\varepsilon,u_0}|^2 u^{\varepsilon,u_0} + \sqrt{\varepsilon}\Gamma, \tag{1.3}$$

where $\varepsilon$ stands for the small noise amplitude, $\Gamma$ is a complex Gaussian space-time white noise and $u_0$ is again the initial datum.

A formal derivation of this model is proposed in the above references. In particular, it is argued that it can be assumed that the damping term is exactly balanced by the amplifiers so that these effects do not appear directly in the model. Only the noise remains. Note also that this equation also appears in the context of anharmonic atomic chains in the presence of thermal fluctuation; see, for example, [5] where timing jitter is also studied. In that second case the derivation can probably be done in a more rigorous way; this will be the object of future work.

If we consider the recently studied Raman coupling to the thermal phonon (see [7, 15, 16, 17, 29]), or four-wave-mixing (see [15, 30]), another quantum noise appears and depends on the pulse intensity. It is modeled as a real multiplicative noise. The first physical derivation of the equation is obtained in [16]. Note that in the case of the Raman amplification, an extra Raman nonlinear response appears in the equation. As in the above physical references, we drop it since it is assumed to have a limited effect on the noise induced timing jitter. The following model is used:

$$i\frac{\partial u^{\varepsilon,u_0}}{\partial t} = \Delta u^{\varepsilon,u_0} + |u^{\varepsilon,u_0}|^2 u^{\varepsilon,u_0} + \sqrt{\varepsilon} u^{\varepsilon,u_0}\Gamma_R. \tag{1.4}$$

Here the noise $\Gamma_R$ is a real Gaussian noise and the product is a Stratonovich product. An important feature of this type of noise is that the mass, given mathematically by the square of the $L^2$ norm, is a conserved quantity. The stochastic NLS equation with real multiplicative noise is also used in the context of crystals; see, for example, [2, 3].

Unlike the deterministic case, an initial soliton profile is progressively distorted due to noise. As a consequence, with a probability that is expected to be small, an error in transmission occurs. It is an important issue to derive theoretical tools to estimate this probability. A first type error occurs when a soliton is emitted and at the other end of the fiber it is not detected. Two phenomena may induce such error. The mass

$$\mathbf{N}(\phi) = \|\phi\|_{\mathrm{L}^2}^2 \qquad \forall \phi \in \mathrm{L}^2,$$

which is an invariant quantity without noise fluctuates when an additive noise is taken into account. Thus, when the noise is additive, the signal may not be detected due to a decrease of the mass. The second source of error is the so-called timing jitter. The arrival time is defined as

$$\mathbf{Y}(\phi) = \int_{\mathbb{R}} x|\phi(x)|^2\, dx \qquad \forall \phi \in \mathrm{L}^2.$$

Without noise, the signal is centered at time $x = 0$ and the arrival time is zero. The noise may change the arrival time and shift the signal outside the measuring window $[-l/2, l/2]$. From these considerations, the problem is reduced to estimate the probability that the mass has decreased significantly or that the arrival time has changed significantly.

Similarly, when no signal is emitted, an additive noise may create from nothing a signal with high enough mass at $T$ and that might be mistaken as a 1. When the noise is multiplicative, because the mass is invariant, we only have to consider the loss of a 1 due to timing jitter.

The aim of the paper is to apply probabilistic tools, more specifically, large deviations estimates to evaluate theoretically the probability of large fluctuations of the mass and arrival time. The large fluctuations events are indexed by $R$ positive, for example, large fluctuations of the arrival time correspond to $\{\mathbf{Y}(u^{\varepsilon,u^0}(T)) \geq R\}$ or $\{\mathbf{Y}(u^{\varepsilon,u_0}(T)) \leq -R\}$. Using large deviation techniques is justified by the standard assumption in the physics literature that the noise is small. We prove that, as usual in that context, the large deviation probabilities are related to an optimal control problem. They are deduced by contraction from a large deviation principle at the level of the paths. Our aim is to give precise upper and lower bounds of these large deviation probabilities.

We get lower bounds by minimizing the rate function over a small set of paths. Namely, we take paths which are modulated solitons, that is, solitons with time varying parameters. Upper bounds are obtained using energy inequalities.

In the physics literature a different method is used. It relies on an adiabatic perturbation theory, (see, e.g., [25, 26]), where the pulse is approximated by a soliton ansatz with finite fluctuating collective variables. In other words, the stochastic NLS equation is replaced by a finite number of coupled stochastic differential equations for the soliton parameters. Thus, the original infinite dimensional problem is reduced to a finite dimensional one for which powerful methods can be used.

It seems very difficult to justify theoretically this method. Our argument is rigorous. Soliton ansatz are also used, but only to provide lower bounds. Surprisingly, the upper bound which is obtained in a totally different way is of the same order as the lower bound with respect to different physically relevant parameters: the length of the fiber $T$, the initial amplitude of the signal $A$ and the parameter $R$ indexing the large fluctuation event. Moreover, our results are comparable to the ones available in the physics literature.

We recover, for instance, the fact that the law of the mass is not Gaussian. Concerning the arrival time, the order in $R$ we obtain proves that the log of the tails are undistinguishable from the log of Gaussian tails. The Gaussianity of the arrival time is a well studied issue in physics. In [24], assuming that the timing jitter is the most troublesome process, an upper limit of the information rate is derived based on the Gaussian assumption and variance computations. In [33] the log of the tails of the amplitude and arrival time are evaluated numerically via an importance sampling Monte Carlo estimator and using an ansatz approximation. It is obtained that the log of tails of the arrival time is the same as the log of a Gaussian tail, while the log of the tail of the amplitude differs significantly from that of Gaussian tails. In [1, 14, 28] it is shown that the arrival time can be considered as Gaussian in the first order only while in [14, 20, 32, 37] justifications for a deviation from Gaussianity when there is filtering or soliton interaction are given.

Assuming that the seemingly Gaussian arrival time is indeed Gaussian, we obtain the same order in $T$ as physicists. In [5, 17, 24] the variance of the arrival time is studied. In [24] it is proved to be of the order of $T^3$ (superdiffusion) and the timing jitter is connected to a shift in the soliton carrier frequency which we exploit in the construction of ansatz. In [17], where both independent complex additive and real multiplicative noises appear in the equation, the contribution of each noise to the variance of the arrival time is also of the order $T^3$.

We are also able to compare the tails of the arrival time to the tails of the mass at the end of the line when the noise is additive. We obtain that the tails of the arrival time are thicker than that of the mass. Thus, timing jitter is the dominant factor as suggested by Gordon and Haus in [24].

Let us also mention that recent articles [14, 19, 31] give approximate PDF of the mass, as well as of the joint law of the mass and arrival time at $T$ with initial datum $\Psi_A^0$. Our results compare to theirs in many ways. In the first article the PDF is obtained using the Fokker–Planck equation and again approximating the pulse by a soliton with finite random modulations evolving according to dynamically coupled SDEs. In [19] the PDF is obtained via a saddle point approximation of a finite dimensional approximation of the infinite dimensional Martin–Siggia–Rose effective action, relying on ansatz. Theory for infinite dimensional effective action is developed in [27], but it has not been used so far for the problem at hand. These infinite dimensional effective actions in physics are intimately related to the rate function of a sample path large deviation principle (LDP). Paths minimizing the action are then called optimal fluctuations or instantons generalizing the quantum mechanics instantons studied in [21] using large deviation techniques.

With our large deviations approach, we study the tails of the CDF and not the bulk of the distribution as with PDFs. The bulk seems less interesting for a study of the rare events causing error in transmission. We obtain accurate rigorous results without using directly the spectral properties of the nonlinear Schrödinger operator. Though applied here to the problem of the error in transmission and for specific and simplified equations, this approach could be used for more general models. Its application to the exit time off neighborhoods of the soliton or randomly modulated soliton for stochastic Korteweg–de Vries equations will be given elsewhere.

**2. Notation and preliminaries.** For $p \geq 1$, $\mathrm{L}^p$ is the classical Lebesgue space of complex valued functions on $\mathbb{R}$ and $\mathrm{W}^{1,p}$ is the associated Sobolev space of $\mathrm{L}^p$ functions with first order derivatives, in the sense of distributions, in $\mathrm{L}^p$. If $I$ is an interval of $\mathbb{R}$, $(E, \|\cdot\|_E)$ a Banach space and $r$ belongs to $[1,\infty]$, then $\mathrm{L}^r(I;E)$ is the space of strongly Lebesgue measurable functions $f$ from $I$ into $E$ (see [18]) such that $t \to \|f(t)\|_E$ is in $\mathrm{L}^r(I)$. The space $\mathrm{L}^2$ with the inner product defined by $(u,v)_{\mathrm{L}^2} = \mathfrak{Re} \int_{\mathbb{R}} u(x)\overline{v}(x)\,dx$ is a Hilbert space. The Sobolev spaces $\mathrm{H}^s$ are the Hilbert spaces of functions of $\mathrm{L}^2$ with partial derivatives up to order $s$ in $\mathrm{L}^2$. When $s$ is fractional it is defined classically via the Fourier transform. When the functions are real valued we specify it, for example, we write $\mathrm{H}^s(\mathbb{R},\mathbb{R})$. The following Hilbert spaces of spatially localized functions

$$\Sigma = \{f \in \mathrm{H}^1 : x \mapsto xf(x) \in \mathrm{L}^2\},$$
$$\Sigma^{1/2} = \{f \in \mathrm{H}^1 : x \mapsto \sqrt{|x|}f(x) \in \mathrm{L}^2\}$$

are also introduced and endowed with the norms

$$\|f\|_{\Sigma}^2 = \|f\|_{\mathrm{H}^1}^2 + \|x \mapsto xf(x)\|_{\mathrm{L}^2}^2,$$
$$\|f\|_{\Sigma^{1/2}}^2 = \|f\|_{\mathrm{H}^1}^2 + \|x \mapsto \sqrt{|x|}f(x)\|_{\mathrm{L}^2}^2.$$

We denote by $\|\Phi\|_{\mathcal{L}_c(A,B)}$ the norm of $\Phi$ as a linear continuous operator from $A$ to $B$, where $A$ and $B$ are normed vector spaces. We recall that $\Phi$ is a Hilbert–Schmidt operator from $H$ to $\tilde{H}$, where $H$ and $\tilde{H}$ are Hilbert spaces, if it is a linear continuous operator such that, given a complete orthonormal system $(e_j^H)_{j=1}^\infty$ of $H$, $\sum_{j=1}^\infty \|\Phi e_j^H\|_{\tilde{H}}^2 < \infty$. We denote by $\mathcal{L}_2(H,\tilde{H})$ the space of Hilbert–Schmidt operators from $H$ to $\tilde{H}$ endowed with the norm

$$\|\Phi\|_{\mathcal{L}_2(H,\tilde{H})} = \mathrm{tr}(\Phi\Phi^*) = \sum_{j=1}^\infty \|\Phi e_j^H\|_{\tilde{H}}^2.$$

We also recall that a cylindrical Wiener process $W_c$ in a Hilbert space $H$ is such that, for any complete orthonormal system $(e_j)_{j=1}^\infty$ of $H$, there exists a sequence of independent Brownian motions $(\beta_j)_{j=1}^\infty$ such that $W_c = \sum_{j=1}^\infty \beta_j e_j$. This sum does not converge in $\mathrm{H}^1$ but in any Hilbert space $U$ such that the embedding $H \subset U$ is Hilbert–Schmidt. The image of the process $W_c$ by a linear mapping $\Phi$ on $H$ is a well defined process in $H$ when the mapping is Hilbert–Schmidt on $H$, that is, $\Phi \in \mathcal{L}_2(H) = \mathcal{L}_2(H,H)$. Then, $W = \Phi W_c$ is such that $W(1)$ is a well defined Gaussian random variable with covariance operator $\Phi\Phi^*$. A detailed presentation of Hilbert space valued Wiener processes, the stochastic integration in that setting and SPDEs is given, for instance, in [8], Chapter 4.

We recall that a rate function $I$ is a lower semicontinuous function and that a good rate function $I$ is a rate function such that, for every positive $c$, $\{x : I(x) \le c\}$ is a compact set.

Let us now recall some mathematical aspects of the stochastic NLS equations. The equations, written as SPDEs in the Itô form, are in the additive case

$$idu^{\varepsilon,u_0} - (\Delta u^{\varepsilon,u_0} + |u^{\varepsilon,u_0}|^2 u^{\varepsilon,u_0})\,dt = \sqrt{\varepsilon}\,dW, \tag{2.1}$$

and in the multiplicative case

$$idu^{\varepsilon,u_0} - (\Delta u^{\varepsilon,u_0} + |u^{\varepsilon,u_0}|^2 u^{\varepsilon,u_0})\,dt = \sqrt{\varepsilon} u^{\varepsilon,u_0} \circ dW. \tag{2.2}$$

The symbol $\circ$ stands for the Stratonovich product. It is convenient to use the Itô product so that we write the equivalent Itô form of the equation

$$idu^{\varepsilon,u_0} - (\Delta u^{\varepsilon,u_0} + |u^{\varepsilon,u_0}|^2 u^{\varepsilon,u_0})\,dt = \sqrt{\varepsilon} u^{\varepsilon,u_0}\,dW - \frac{i}{2}\varepsilon F_\Phi u^{\varepsilon,u_0}, \tag{2.3}$$

where, given $(e_j)_{j=1}^\infty$ an orthonormal basis of $\mathrm{L}^2$, $F_\Phi(x) = \sum_{j=1}^\infty (\Phi e_j)^2(x)$. The term $(\varepsilon/2)F_\Phi(x)$ is the Itô correction necessary to transform the Stratonovich product into a Itô one. Note that $F_\Phi$ does not depend on the basis.

As mentioned earlier, in the case of equation (2.3) (see [10]), the mass

$$\mathbf{N}(u^{\varepsilon,u_0}(t)) = \|u^{\varepsilon,u_0}(t)\|_{\mathrm{L}^2}^2, \qquad t > 0,$$

is a conserved quantity. Precise assumptions on $\Phi$ such that $W = \Phi W_c$ are made below. These equations are supplemented with an initial datum

$$u^{\varepsilon,u_0}(0) = u_0.$$

In this paper we consider initial data in $\Sigma \subset \mathrm{H}^1$ and work with the solution constructed in [10]. Since we work with a subcritical nonlinearity, we could also consider solutions in $\mathrm{L}^2$ with initial data in $\mathrm{L}^2$. However, the $\mathrm{H}^1$-setting is preferred in order to be able to consider the spaces $\Sigma$ and $\Sigma^{1/2}$ and study the arrival time of the pulse

$$\mathbf{Y}(u^{\varepsilon,u_0}(t)) = \int_{\mathbb{R}} x |u^{\varepsilon,u_0}(t,x)|^2\,dx, \qquad t \geq 0,$$

defined when $u^{\varepsilon,u_0}(t)$ belongs to $\Sigma^{1/2}$.

We are concerned by weak solutions or, equivalently, by mild solutions which, in the additive case, satisfy

$$\begin{aligned} u^{\varepsilon,u_0}(t) = U(t)u_0 - i\int_0^t U(t-s)(|u^{\varepsilon,u_0}(s)|^2 u^{\varepsilon,u_0}(s))\,ds \\ -i\sqrt{\varepsilon}\int_0^t U(t-s)\,dW(s), \end{aligned} \tag{2.4}$$

where $(U(t))_{t\in\mathbb{R}}$ stands for the Schrödinger group, $U(t)=e^{-it\Delta}, t\in\mathbb{R}$. The last term is called the stochastic convolution. In the multiplicative case, the mild equation is

$$\begin{aligned} u^{\varepsilon,u_0}(t) = U(t)u_0 - i\int_0^t U(t-s)(|u^{\varepsilon,u_0}(s)|^2 u^{\varepsilon,u_0}(s))\,ds \\ -i\sqrt{\varepsilon}\int_0^t U(t-s)u^{\varepsilon,u_0}(s)\,dW(s) \\ +(\varepsilon/2)\int_0^t U(t-s)F_\Phi u^{\varepsilon,u_0}(s)\,ds, \end{aligned} \tag{2.5}$$

where the stochastic integral is an Itô integral.

The noise is the time derivative in the sense of distributions of the Wiener process $W$. It corresponds to a white noise in time. A space-time white noise would correspond to $\Phi$ equal to the identity. It is, however, the noise mainly considered in optics. We cannot handle such rough noises and make the assumption that the noises are colored in space in order to obtain well-posedness. The basic limitation is that, unlike semi-groups like the Heat semi-group, the Schrödinger group is an isometry and does not allow smoothing in the Sobolev spaces based on $\mathrm{L}^2$. For instance, in the additive case, it can be seen that the stochastic convolution is a well defined process with paths in $\mathrm{L}^2$ if and only if $\Phi$ is a Hilbert–Schmidt operator on $\mathrm{L}^2$. In that case, however, we will, for computational issues, consider sequences of noises that mimic the white noise in the limit. This statement will be made more precise.

In fact, we make even stronger assumptions. In the additive case we assume that $W$ is a Wiener process in $\Sigma$, in other words, we require that $\Phi\in\mathcal{L}_2(L^2,\Sigma)$. In the multiplicative case, it is imposed that $W$ is a Wiener process in $\mathrm{H}^s(\mathbb{R},\mathbb{R})$ where $s$ satisfies $s>3/2$. It allows to consider paths in $\Sigma^{1/2}$.

We know that the Cauchy problem is globally well posed in $\mathrm{H}^1$; see [10] for a general discussion on the local well posedness and the global existence for more general nonlinearities and dimensions. Note that in the deterministic case, the NLS equation considered here is integrable thanks to the inverse scattering method. We do not use these techniques in this article. Results on the influence of the noise on the blow-up time for more general nonlinearities and dimensions are given in [11, 12]. In [4, 13] the ideal white noise and results on the influence of a noise on the blow-up are studied numerically.

Sample path LDPs for stochastic NLS equations are proved in [22, 23]. These LDPs are stated in the topology of $\mathrm{C}([0,T];\mathrm{H}^1)$ for $T>0$ and do not allow to treat the arrival time of the solution. We shall generalize these and give LDPs in $\mathrm{C}([0,T];\Sigma^{1/2})$. The rate function of the LDP in the additive case is defined in terms of the mild solution of the control problem

$$\begin{cases} i\dfrac{du}{dt} = \Delta u + |u|^2 u + \Phi h, \\ u(0) = u_0 \in \Sigma \text{ and } h \in \mathrm{L}^2(0,T;\mathrm{L}^2). \end{cases} \tag{2.6}$$

We denote the solution by $u = \mathbf{S}^{a,u_0}(h)$. The mapping $h \to \mathbf{S}^{a,u_0}(h)$ is called the control map and (2.6) the control equation.

In the multiplicative case, the control equation is

$$i\frac{du}{dt} = \Delta u + |u|^2 u + u\Phi h, \tag{2.7}$$

whose mild solution is denoted by $u = \mathbf{S}^{m,u_0}(h)$. The mapping $\mathbf{S}^{m,u_0}$ is again the control map and (2.7) the control equation.

In this article, when describing properties which hold both in the additive and multiplicative cases, we use the symbol $\mathbf{S}(u_0, h)$ to denote either $\mathbf{S}^{a,u_0}(h)$ or $\mathbf{S}^{m,u_0}(h)$.

Let us now state the sample path LDPs. As already mentioned, these are slight generalizations of the LDPs given in [22, 23]. For the reader's convenience, we give the proof for the case of an additive noise in Section 5. The case of a multiplicative noise is more involved, but does not present new difficulties compared to the proof given in [23]. In order to keep the length of the article reasonable, we only give the new ingredients necessary to adapt the proof.

Theorem 2.1. *Assume that $\Phi$ belongs to $\mathcal{L}_2(\mathrm{L}^2, \Sigma)$ in the additive case and $\Phi \in \mathcal{L}_2(\mathrm{L}^2, \mathrm{H}^s(\mathbb{R}, \mathbb{R}))$ with $s > 3/2$ in the multiplicative case. Assume also that the initial datum $u_0$ is in $\Sigma$. Then the solutions of the stochastic nonlinear Schrödinger equations (2.4) and (2.5) are almost surely in $\mathrm{C}([0,T];\Sigma^{1/2})$. Moreover, they define $\mathrm{C}([0,T];\Sigma^{1/2})$ random variables and their laws $(\mu^{u^{\varepsilon,u_0}})_{\varepsilon>0}$ satisfy a LDP of speed $\varepsilon$ and good rate function*

$$I^{u_0}(w) = \tfrac{1}{2} \inf_{h \in \mathrm{L}^2(0,T;\mathrm{L}^2): w = \mathbf{S}(u_0,h)} \|h\|^2_{\mathrm{L}^2(0,T;\mathrm{L}^2)},$$

*where $\mathbf{S}(u_0,\cdot) = \mathbf{S}^{a,u_0}(\cdot)$ in the additive case and $\mathbf{S}(u_0,\cdot) = \mathbf{S}^{m,u_0}(\cdot)$ in the multiplicative case, and with the convention that $\inf \varnothing = \infty$. It means that, for every Borel set $B$ of $\mathrm{C}([0,T];\Sigma^{1/2})$, we have the lower bound*

$$-\inf_{w \in \overset{\circ}{B}} I^{u_0}(w) \le \varliminf_{\varepsilon \to 0} \varepsilon \log \mathbb{P}(u^{\varepsilon,u_0} \in B)$$

*and the upper bound*

$$\varlimsup_{\varepsilon \to 0} \varepsilon \log \mathbb{P}(u^{\varepsilon,u_0} \in B) \le -\inf_{w \in \overline{B}} I^{u_0}(w).$$

These sample path LDPs allow, for example, to evaluate the probability that, originated from a soliton profile

$$\Psi_A^0(x) = \sqrt{2}A\,\mathrm{sech}(Ax),$$

the random solution be significantly different from the deterministic soliton solution

$$\Psi_A(t,x) = \Psi_A^0(x)\exp(-iA^2t).$$

Indeed, for $\delta$ and $\eta$ positive and $\varepsilon$ small enough, the LDP implies that

$$\exp(-C_1/\varepsilon) \le \mathbb{P}(\|u^{\varepsilon,\Psi_A^0} - \Psi_A\|_{\mathrm{C}([0,T];\Sigma^{1/2})} > \delta) \le \exp(-C_2/\varepsilon),$$

where

$$C_1 = \inf_{w:\|w-\Psi_A\|_{\mathrm{C}([0,T];\Sigma^{1/2})}>\delta} I^{\Psi_A^0}(w) + \eta$$

and

$$C_2 = \inf_{w:\|w-\Psi_A\|_{\mathrm{C}([0,T];\Sigma^{1/2})}\ge\delta} I^{\Psi_A^0}(w) - \eta.$$

Recall that, since the rate function is a good rate function, if $B$ is a closed set and $\inf_{w\in B} I^{\Psi_A^0}(w) < \infty$, then there is an $f$ in $B$, optimal fluctuation, such that $I^{\Psi_A^0}(f) = \inf_{w\in B} I^{\Psi_A^0}(w)$. Thus, if $B$ does not contain the deterministic solution, then necessarily $\inf_{w\in B} I^{\Psi_A^0}(w) > 0$. Consequently, $\eta$ may be chosen such that $C_2$ is positive and the above probability of a deviation from the deterministic path is exponentially small in the small $\varepsilon$ limit.

In this article we are interested in estimating the probability of particular deviations from the deterministic paths. Namely, we wish to study how the mass and the arrival time of a solution at coordinate $T$ deviate from their value in the "frozen" deterministic system (i.e., when $\varepsilon = 0$). In the absence of noise, the mass is a conserved quantity and for initial data being either 0 or $\Psi_A^0$ the arrival time remains equal to zero.

We know from [22] that we may use the contraction principle to deduce from LDP for the paths a LDP for the mass at $T$ and obtain a LDP with speed $\varepsilon$ and good rate function for an initial datum $u_0$

$$I_N^{u_0}(m) = \frac{1}{2} \inf_{h\in\mathrm{L}^2(0,T;\mathrm{L}^2):\mathbf{N}(\mathbf{S}^{a,u_0}(h)(T))=m} \{\|h\|^2_{\mathrm{L}^2(0,T;\mathrm{L}^2)}\}.$$

In the case of a multiplicative noise, the mass is a conserved quantity. Thus, in this case, the mass cannot deviate from the constant value corresponding to that of the initial datum.

Similarly, the mapping $\mathbf{Y}$ is continuous from $\Sigma^{1/2}$ into $\mathbb{R}$. We may thus define by direct image the measures $(\mu^{\mathbf{Y}(u^{\varepsilon,u_0}(T))})_{\varepsilon>0}$ for an initial datum $u_0$ in $\Sigma$. We obtain by contraction that they satisfy a LDP of speed $\varepsilon$ and good rate function

$$I_Y^{u_0}(y) = \frac{1}{2} \inf_{h \in \mathrm{L}^2(0,T;\mathrm{L}^2): \mathbf{Y}(\mathbf{S}(u_0,h)(T))=y} \{\|h\|^2_{\mathrm{L}^2(0,T;\mathrm{L}^2)}\},$$

the control map $\mathbf{S}$ is either that of the additive or multiplicative case.

Let us briefly explain our strategy to estimate the probability of some event. Let us consider, for instance, the event $D_\varepsilon = \{\mathbf{Y}(u^{\varepsilon,0}(T)) \in [a,b]\}$, where $[a,b]$ is an interval which does not contain 0. We use the LDP to obtain

$$\begin{aligned} -\inf_{y\in(a,b)} I_Y^0(y) &\le \varliminf_{\varepsilon\to 0} \varepsilon \log \mathbb{P}(D_\varepsilon) \\ &\le \varlimsup_{\varepsilon\to 0} \varepsilon \log \mathbb{P}(D_\varepsilon) \le -\inf_{y\in[a,b]} I_Y^0(y). \end{aligned} \tag{2.8}$$

To approximate from above the upper bound, we use energy type inequalities. These give a lower bound on the minimum $\mathrm{L}^2$ norm of the control $h$ required to change the deterministic behavior and have the arrival time in $[a,b]$ at time $T$. Namely, we obtain a positive constant $c$ such that

$$\text{if } \mathbf{Y}(S(u_0,h)(T)) \in [a,b] \qquad \text{then } \tfrac{1}{2}\|h\|^2_{\mathrm{L}^2(0,T;\mathrm{L}^2)} \ge c.$$

This clearly implies

$$\varlimsup_{\varepsilon\to 0} \varepsilon \log \mathbb{P}(D_\varepsilon) \le -c.$$

The second step is to find a particular function $h$ such that $\mathbf{Y}(S(u_0,h)(T)) \in (a,b)$ and $\tilde{c} = (1/2)\|h\|^2_{\mathrm{L}^2(0,T;\mathrm{L}^2)}$ is as small as possible. Then

$$-\tilde{c} \le \varliminf_{\varepsilon\to 0} \varepsilon \log \mathbb{P}(D_\varepsilon).$$

In this second step, we minimize on a smaller set of controls which gives rise to a problem a from the calculus of variations.

The difficulty is to have sufficiently sharp energy estimates and to find a good solution to the control problem so that $c$ and $\tilde{c}$ are as close as possible. We see below that we are able to do so in some interesting situations and derive good estimates on such probabilities.

Note that proceeding as in [22] for the mass, we may prove in the additive case that $\inf_{y\in J} I_Y^{u_0}(y) < \infty$ for every nonempty interval $J$ and any $u_0$ provided the range of $\Phi$ is dense. Indeed, for every real number $a$, a solution of the form $u(t,x) = (1+atx)u_0$ satisfies $\mathbf{Y}(u(T)) = aT\pi^2/3$. Plugging this solution into equation (2.6), we find a control such that the solution reaches any interval at time $T$. Using the continuity of $h \mapsto \mathbf{Y}(\mathbf{S}^{a,u_0}(h)(T))$ from $\mathrm{L}^2(0,T;\mathrm{L}^2)$ into $\mathbb{R}$ and the density of the range of $\Phi$, we obtain $\inf_{y\in J} I_Y^{u_0}(y) < \infty$. This shows that the lower bound is nontrivial. With arguments given previously, we know that the upper bound is nontrivial as well.

REMARK 2.2. Using similar arguments as in [22], we can prove that, for every positive $R$ besides an at most countable set of points, we can replace $\underline{\lim}$ and $\overline{\lim}$ by lim in the LDP and obtain

$$\begin{aligned}
&\lim_{\varepsilon\to 0} \varepsilon \log \mathbb{P}(\mathbf{Y}(u^{\varepsilon,u_0}(T)) \geq R) \\
&\qquad = -(1/(2\varepsilon)) \inf_{h\in \mathrm{L}^2(0,T;\mathrm{L}^2):\mathbf{Y}(\mathbf{S}(u_0,h)(T))\geq R} \{\|h\|^2_{\mathrm{L}^2(0,T;\mathrm{L}^2)}\}, \\
&\lim_{\varepsilon\to 0} \varepsilon \log \mathbb{P}(\mathbf{Y}(u^{\varepsilon,u_0}(T)) \leq -R) \\
&\qquad = -(1/(2\varepsilon)) \inf_{h\in \mathrm{L}^2(0,T;\mathrm{L}^2):\mathbf{Y}(\mathbf{S}(u_0,h)(T))\leq -R} \{\|h\|^2_{\mathrm{L}^2(0,T;\mathrm{L}^2)}\}.
\end{aligned}$$

This uses the fact that a monotone and bounded function is continuous almost everywhere.

We end this section with some remarks which will be useful in the development of our method when we consider the arrival time of the solution. Let us consider the case when the initial datum is $\Psi_A^0$. The probability of tail events of the arrival time are related to the behavior of $\mathbf{Y}(S(\Psi_A^0, h))$.

An equation for the motion of the arrival time is given in [38] in the case of an external potential. The first step consists in multiplying the control PDE by $-ix\overline{u}$, taking the real part, and integrating by part the term involving the Laplace operator. We then obtain for the control PDE associated to the multiplicative case

$$\frac{d}{dt}\mathbf{Y}(\mathbf{S}^{m,\Psi_A^0}(h)(t)) = 2\mathfrak{Re}\left(i\int_{\mathbb{R}} \overline{\mathbf{S}^{m,\Psi_A^0}(h)(t,x)}\, \partial_x \mathbf{S}^{m,\Psi_A^0}(h)(t,x)\, dx\right), \tag{2.9}$$

while in the additive case we obtain

$$\begin{aligned}
\frac{d}{dt}\mathbf{Y}(\mathbf{S}^{a,\Psi_A^0}(h)(t)) &= 2\mathfrak{Re}\left(i\int_{\mathbb{R}} \overline{\mathbf{S}^{a,\Psi_A^0}(h)(t,x)} \partial_x \mathbf{S}^{a,\Psi_A^0}(h)(t,x)\, dx\right) \\
&\quad - 2\mathfrak{Re}\left(i\int_{\mathbb{R}} x\overline{\mathbf{S}^{a,\Psi_A^0}(h)(t,x)}(\Phi h)(t,x)\, dx\right).
\end{aligned} \tag{2.10}$$

Below, we use the notation

$$\mathbf{P}(u) = 2\mathfrak{Re}\left(i\int_{\mathbb{R}} \overline{u}(x)\partial_x u(x)\, dx\right), \qquad u \in \mathrm{H}^1.$$

As a consequence of (2.9), we see that, in the multiplicative case, the arrival time of the solution of the control problem cannot move unless its phase depends on the space variable. For instance, if the control is chosen so that the solution $\mathbf{S}(\Psi_A^0, h)$ is a modulated soliton of type (1.2) with varying amplitude and group velocity,

$$\begin{aligned}
\mathbf{S}(\Psi_A^0, h)(t) &= \sqrt{2}A(t)\,\mathrm{sech}(A(t)(x-x_0) + 2A(t)V(t)t) \\
&\qquad \times \exp(-i(A(t)^2 - V(t)^2)t + iV(t)(x-x_0) + i\theta_0),
\end{aligned}$$

we have the well-known identity

$$\frac{d}{dt}\mathbf{Y}(\mathbf{S}(\Psi_A^0,h)(t)) = -2V(t)\mathbf{N}(\mathbf{S}^{m,\Psi_A^0}(h)(t)) = -8V(t)A(t).$$

It will be convenient to choose controlled solutions of the form above. Since the initial datum is $\Psi_A^0$, we necessarily have $V(0)=0$, hence, $V$ cannot be chosen constant, otherwise the arrival time does not change. We will see that it is sufficient to have a constant, amplitude $A$ in order to get sharp bounds. Thus, we will use modulated solitons as solutions of the control problem with constant amplitude when studying the motion of the arrival time.

The first idea to find a control giving a solution whose arrival time or mass verify some desired property is to take the above modulated soliton and plug it into the control equation. This gives an explicit form of the control in terms of the various parameters. Then, we compute the space-time $\mathrm{L}^2$ norm of this control. We obtain a function of the parameters which we can try to minimize thanks to the calculus of variations. This approach is not easy to perform, the function to minimize has a complicated form and is often singular. However, we are not interested in finding the exact extremal of the function and we use this method in an heuristic way. This allows us to guess good controls with low enough space-time $\mathrm{L}^2$ norm compared to the upper bound and such that the controlled path is a modulated soliton satisfying the desired constraints. We will see that this method is successful.

Let us consider the following controlled nonlinear Schrödinger equation:

$$i\frac{du}{dt} = \Delta u + |u|^2u + \lambda(t)xu \tag{2.11}$$

with initial datum $\Psi_A^0$. The function $\lambda$ is taken in $\mathrm{L}^1(0,T;\mathbb{R})$. This corresponds to the multiplicative control equation with $\Phi h = \lambda(t)x$ or to the additive one with $\Phi h = \lambda(t)xu$. We use well-known transformations to compute explicitly the solution of (2.11) which we denote by $\Psi_{A,\lambda}$. We first may check that the functions $v_1$ and $v_2$ defined by $v_1(t,x) = \exp(i(\int_0^t \lambda(s)\,ds)x)u(t,x)$ and $v_2(t,x) = \exp(-i\int_0^t(\int_0^s \lambda(\tau)\,d\tau)^2\,ds)v_1(t,x)$ (gauge transform) satisfy the PDEs

$$i\frac{\partial v_1}{\partial t} = \frac{\partial^2 v_1}{\partial x^2} + |v_1|^2v_1 - \left(\int_0^t \lambda(s)\,ds\right)^2 v_1 - 2i\left(\int_0^t \lambda(s)\,ds\right)\frac{\partial v_1}{\partial x}$$

and

$$i\left(\frac{\partial v_2}{\partial t} + 2\left(\int_0^t \lambda(s)\,ds\right)\frac{\partial v_2}{\partial x}\right) = \frac{\partial^2 v_2}{\partial x^2} + |v_2|^2 v_2$$

with initial datum $\Psi_A^0$. We conclude using the methods of characteristics that $v_3$ defined by

$$v_3(t,x) = v_2\left(t, x + 2\int_0^t\int_0^s \lambda(u)\,du\,ds\right)$$

is a solution of the usual NLS equation with initial datum $\Psi_A^0$. Thus, we obtain that $v_3(t,x)=\Psi_A(t,x)$ and that the solution of the Cauchy problem associated to (2.11) is

$$\begin{aligned}\Psi_{A,\lambda}(t,x) = \sqrt{2}A\,\mathrm{sech}\Big(A\Big(x-2\int_0^t\int_0^s\lambda(\tau)\,d\tau\,ds\Big)\Big)\\ \times\exp\Big[-iA^2t+i\int_0^t\Big(\int_0^s\lambda(\tau)\,d\tau\Big)^2ds\\ -ix\int_0^t\lambda(s)\,ds+2i\Big(\int_0^t\lambda(s)\,ds\Big)\Big(\int_0^t\int_0^s\lambda(\tau)\,d\tau\,ds\Big)\Big].\end{aligned}$$

We obtain a modulated soliton with group velocity given by $V(t)=\int_0^t\lambda(s)\,ds$. In the additive case, it is possible to obtain a control such that the solution has the same arrival time and group velocity and such that the space-time $L^2$ norm of the control is simpler to compute. It is obtained thanks to the observation that using the gauge transform the solution of the Cauchy problem

$$\begin{cases} i\dfrac{dv}{dt}=\Delta v+|v|^2v+\lambda(t)\Big(x-2\displaystyle\int_0^t\int_0^s\lambda(\tau)\,d\tau\,ds\Big)v,\\ v(0)=\Psi_A^0,\end{cases}\tag{2.12}$$

is given by

$$\tilde{\Psi}_{A,\lambda}(t,x)=\exp\Big(2i\int_0^t\lambda(s)\int_0^s\int_0^\tau\lambda(\sigma)\,d\sigma\,d\tau\,ds\Big)\Psi_{A,\lambda}(t,x).$$

Remark 2.3. For the controls chosen above, relation (2.9) holds also in the additive case. Thus, the second term in (2.10) which, at first glance, could be useful to act on the arrival time is in fact useless.

Also, it could be thought that the choice of more complicated group velocities could be useful. We have tried to consider a space dependent group velocity, but the calculus of variations approach indicates that optimality is reached when it does not depend on space.

**3. Tails of the mass and arrival time with additive noise.** In the case of an additive noise, both the mass and arrival time may deviate from the deterministic behavior and result in error in transmission.

We study tails and thus the probability of a deviation from the mean. The constant $R$ below quantify this deviation. We are not really interested in large $R$. In the case of the mass, for example, interesting cases are when $R$ lies in $(0,4)$. But, since $\varepsilon$ goes to zero and the factor in the exponential is of the order of $1/\varepsilon$ while $R$ is of order 1, it results in very unlikely events. These significant excursions of the mass and arrival time are exactly large deviation events.

Moreover, another parameter is particularly interesting. It is $T$ the length of the fiber optical line. It is assumed to be large. For example, we could think of a fiber optical line between Europe and America.

We first recall the results obtained in [22] for the tails of mass of the pulse at the end of the line. We repeat the proofs for the reader's convenience. The aim is to compare these tails with the tails of the arrival time obtained thereafter. We show that indeed the timing jitter is the dominant effect in the error in transmission when the noise is additive. The initial datum may be $u_0 = 0$ or $u_0 = \Psi$, where $\Psi(x) = \sqrt{2}\,\mathrm{sech}(x)$. We could consider a soliton profile with any amplitude $A$ as well but, for simplicity, we consider the case $A = 1$. However, we consider below the parameter $A$ for the timing jitter in order to compare with results from physics. Indeed, this dependance is made explicit in the case of the timing jitter in the physics literature. Let us begin with upper bounds of the tails. As already mentioned, they are obtained thanks to energy estimates. For the second bound, we consider the case of the emission of a signal. In that case only a decrease of the mass is troublesome and causes error in transmission. Thus, the bound given only accounts for a significant decrease of the mass.

Proposition 3.1. *For every positive $T$ and $R$ [$R$ in $(0,4)$ for the second inequality] and every operator $\Phi$ in $\mathcal{L}_2(\mathrm{L}^2, \mathrm{H}^1)$, the following inequalities hold:*

$$\overline{\lim_{\varepsilon\to 0}}\, \varepsilon \log \mathbb{P}(\mathbf{N}(u^{\varepsilon,0}(T)) \geq R) \leq -R/(8T\|\Phi\|^2_{\mathcal{L}_c(\mathrm{L}^2,\mathrm{L}^2)}),$$

$$\overline{\lim_{\varepsilon\to 0}}\, \varepsilon \log \mathbb{P}(\mathbf{N}(u^{\varepsilon,\Psi}(T)) - 4 \leq -R) \leq -R^2/(8T\|\Phi\|^2_{\mathcal{L}_c(\mathrm{L}^2,\mathrm{L}^2)}(4+R)).$$

Proof. Multiplying by $-i\overline{u}$ the equation

$$i\frac{du}{dt} - \Delta u - \lambda |u|^2 u = \Phi h,$$

integrating over time and space and taking the real part gives, for $t \in [0,T]$,

$$\|u(t)\|^2_{\mathrm{L}^2} - \|u_0\|^2_{\mathrm{L}^2} = 2\mathfrak{Re}\left(-i\int_0^t \int_{\mathbb{R}} ((\Phi h)(s,x)\overline{u(s,x)})\,dx\,ds\right) \tag{3.1}$$

and by the Cauchy–Schwarz inequality,

$$\|u(t)\|^2_{\mathrm{L}^2} - \|u_0\|^2_{\mathrm{L}^2} \leq 2\|\Phi\|_{\mathcal{L}_c(\mathrm{L}^2,\mathrm{L}^2)}\|h\|_{L^2(0,T;\mathrm{L}^2)}\|u\|_{L^2(0,T;\mathrm{L}^2)}. \tag{3.2}$$

We integrate once more with respect to $t \in [0,T]$ and obtain

$$\|u\|^2_{L^2(0,T;\mathrm{L}^2)} - T\|u_0\|^2_{\mathrm{L}^2} \leq 2T\|\Phi\|_{\mathcal{L}_c(\mathrm{L}^2,\mathrm{L}^2)}\|h\|_{L^2(0,T;\mathrm{L}^2)}\|u\|_{L^2(0,T;\mathrm{L}^2)}. \tag{3.3}$$

For the first inequality, $u_0 = 0$ and $u = \mathbf{S}^{a,0}(h)$. By (3.3),

$$\|\mathbf{S}^{a,0}(h)\|_{L^2(0,T;\mathrm{L}^2)}^2 \le 2T\|\Phi\|_{\mathcal{L}_c(\mathrm{L}^2,\mathrm{L}^2)}\|h\|_{L^2(0,T;\mathrm{L}^2)}.$$

Then, taking $t = T$ in (3.2), we deduce

$$\|\mathbf{S}^{a,0}(h)(T)\|_{\mathrm{L}^2}^2 \le 4T\|\Phi\|_{\mathcal{L}_c(\mathrm{L}^2,\mathrm{L}^2)}^2\|h\|_{L^2(0,T;\mathrm{L}^2)}^2.$$

Thus, if $N(\mathbf{S}^{a,0}(h)(T)) = \|\mathbf{S}^{a,0}(h)(T)\|_{\mathrm{L}^2}^2 = m$, then

$$\|h\|_{\mathrm{L}^2(0,T;\mathrm{L}^2)}^2 \ge \frac{m}{4T\|\Phi\|_{\mathcal{L}_c(\mathrm{L}^2,\mathrm{L}^2)}^2}.$$

It follows

$$\begin{aligned} I_N^0(m) &= (1/2)\inf_{h\in\mathrm{L}^2(0,T;\mathrm{L}^2):\mathbf{N}(\mathbf{S}^{a,0}(h)(T))=m}\{\|h\|_{\mathrm{L}^2(0,T;\mathrm{L}^2)}^2\} \\ &\ge m/(8T\|\Phi\|_{\mathcal{L}_c(\mathrm{L}^2,\mathrm{L}^2)}^2). \end{aligned}$$

Now, by the LDP for the mass, we have

$$\overline{\lim_{\varepsilon\to 0}}\,\varepsilon\log\mathbb{P}(\mathbf{N}(u^{\varepsilon,0}(T)) \ge R) \le -\inf_{m\in[R,\infty]} I_N^{u_0}(m)$$

and the result follows.

For the second inequality, $u_0 = \Psi$ and $u = \mathbf{S}^{a,\Psi}(h)$. Since $\|\Psi\|_{\mathrm{L}^2}^2 = 4$, (3.3) rewrites

$$\|u\|_{L^2(0,T;\mathrm{L}^2)}^2 - 2T\|\Phi\|_{\mathcal{L}_c(\mathrm{L}^2,\mathrm{L}^2)}\|h\|_{L^2(0,T;\mathrm{L}^2)}\|u\|_{L^2(0,T;\mathrm{L}^2)} - 4T \le 0.$$

Therefore,

$$\begin{aligned} \|u\|_{L^2(0,T;\mathrm{L}^2)} \le\, & T\|\Phi\|_{\mathcal{L}_c(\mathrm{L}^2,\mathrm{L}^2)}\|h\|_{\mathrm{L}^2(0,T;\mathrm{L}^2)} \\ & \times\left(1+\sqrt{1+\frac{4}{T\|\Phi\|_{\mathcal{L}_c(\mathrm{L}^2,\mathrm{L}^2)}^2\|h\|_{\mathrm{L}^2(0,T;\mathrm{L}^2)}^2}}\right). \end{aligned}$$

By (3.2) with $t = T$ we deduce that if $N(\mathbf{S}^{a,0}(h)(T)) = \|\mathbf{S}^{a,0}(h)(T)\|_{\mathrm{L}^2}^2 \le 4 - R$, then

$$R \le 2T\|\Phi\|_{\mathcal{L}_c(\mathrm{L}^2,\mathrm{L}^2)}^2\|h\|_{L^2(0,T;\mathrm{L}^2)}^2\left(1+\sqrt{1+\frac{4}{T\|\Phi\|_{\mathcal{L}_c(\mathrm{L}^2,\mathrm{L}^2)}^2\|h\|_{\mathrm{L}^2(0,T;\mathrm{L}^2)}^2}}\right).$$

We finally obtain

$$\|h\|_{\mathrm{L}^2(0,T;\mathrm{L}^2)}^2 \ge \frac{R^2}{4T\|\Phi\|_{\mathcal{L}_c(\mathrm{L}^2,\mathrm{L}^2)}^2(4+R)}.$$

The upper bound follows. □

Let us now consider lower bounds. As explained above, our method is to find solutions of the control problem with mass at coordinate $T$ satisfying constraints and such that the $L^2$ norm of the control is as small as possible. We find these controlled solutions in the form of modulated solitons. We have found that it is sufficient that only the amplitude varies. We take the solution of (2.6) of the form

$$\sqrt{2}A(t)\exp\left(-i\int_0^t A^2(s)\,ds\right)\operatorname{sech}(A(t)x). \tag{3.4}$$

This is associated to the control

$$\begin{aligned}\Phi h_A(t,x) &= i(A'/A)(t)\Psi_A(t,x)\\ &\quad - i\sqrt{2}A'(t)\exp\left(-i\int_0^t A^2(s)\,ds\right)A(t)x(\sinh/\cosh^2)(A(t)x).\end{aligned}$$

Unfortunately, the right-hand side is in general not in the range of $\Phi$. Moreover, unless we make artificial assumptions on $\Phi$, it is not possible to get information on the norm of $h$. In our result below, we proceed by approximation and consider a sequence of operators $\Phi_n$ approximating the identity on a sufficiently large set containing the controls.

Let us assume for the moment that we can consider the space-time white noise. Then, the mass of the solution (3.4) at time $T$ is equal to $4A(T)$ and the $L^2$ norm of the control is given by

$$\|h_A\|^2_{L^2(0,T;L^2)} = \frac{1}{9}(12+\pi^2)\int_0^T \frac{(A'(t))^2}{A(t)}\,dt. \tag{3.5}$$

The Euler–Lagrange equation associated to the problem of minimizing this quantity is

$$2\frac{A''}{A} = \left(\frac{A'}{A}\right)^2.$$

Multiplying this identity by $A^2$ and differentiating, we obtain $A''' = 0$, so that $A$ is a second degree polynomial and it is easy to see that it has to be of the form $A_0(t) = \alpha(t-\beta)^2$.

For the problem of the zero initial boundary condition, we have $A_0(0) = 0$ and $4A_0(T) > R$. Hence, we deduce the candidates $A_0(t) = \tilde{R}(\frac{t}{2T})^2$ for $\tilde{R} > R$ arbitrary and $R$ defined as for the upper bounds. Plugging such a function into (3.5) gives

$$\|h_{A_0}\|^2_{L^2(0,T;L^2)} = \frac{1}{9}(12+\pi^2)\frac{\tilde{R}}{T}. \tag{3.6}$$

This would give immediately a lower bound if $\Phi$ were the identity. Since we cannot treat this case, we proceed by approximation. The assumption we make on the covariance involves the following sets of time dependent functions. We first introduce, for $D\subset[R,R+1]$,

$$\mathcal{A}_D^1=\{A:[0,T]\to\mathbb{R},\ \text{there exists } \tilde{R}\in D \text{ such that } A(t)=\tilde{R}(t/(2T))^2\}.$$

The functions in $\mathcal{A}_D^1$ are the varying amplitude of the solutions corresponding to controls in the set

$$\mathcal{C}_D^1=\Big\{h\in \mathrm{L}^2(0,T;\mathrm{L}^2),\ \text{there exists } A\in\mathcal{A}_D^1$$
$$h(t,x)=i(A'/A)(t)\Psi_A(t,x)$$
$$-i\sqrt{2}A'(t)\exp\Big(-i\int_0^t A^2(s)\,ds\Big)A(t)x(\sinh/\cosh^2)(A(t)x)\Big\}.$$

For the case of a soliton profile as initial data, a similar argument leads us to define

$$\mathcal{A}_D^2=\{A:[0,T]\to\mathbb{R},\ \text{there exists } \tilde{R}\in D \text{ such that}$$
$$A(t)=(8-\tilde{R}-4\sqrt{4-\tilde{R}})(t/(2T))^2+(-4+2\sqrt{4-\tilde{R}})(t/(2T))+1\}.$$

The set of controls $\mathcal{C}_D^2$ is defined as above by replacing $\mathcal{A}_D^1$ by $\mathcal{A}_D^2$.

We have the following proposition from [22]. The assumptions can easily be fulfilled. They are such that the noise is as close as possible to the space-time white noise considered in physics that we are not able to treat mathematically.

Proposition 3.2. *Let $T$ and $R$ be positive numbers [$R$ in $(0,4)$ for the second inequality], take $D$ dense in $(R,R+1)$ and a sequence of operators $(\Phi_n)_{n\in\mathbb{N}}$ in $\mathcal{L}_2(\mathrm{L}^2,\mathrm{L}^2)$ such that for every $h\in\mathcal{C}_D^1$ we have $\Phi_n h$ converges to $h$ in $\mathrm{L}^1(0,T;\mathrm{L}^2)$. Then we obtain*

$$\underline{\lim}_{n\to\infty,\varepsilon\to0}\ \varepsilon\log\mathbb{P}(\mathbf{N}(u^{\varepsilon,0,n}(T))\ge R)\ge -R(12+\pi^2)/(18T).$$

*Replacing in the above $\mathcal{C}_D^1$ by $\mathcal{C}_D^2$ we obtain*

$$\underline{\lim}_{n\to\infty,\varepsilon\to0}\ \varepsilon\log\mathbb{P}(\mathbf{N}(u^{\varepsilon,\Psi,n}(T))-4\le -R)$$
$$\ge -2(8-R-4\sqrt{4-R})(12+\pi^2)/(36T).$$

*The exponent $n$ is there to recall that $\Phi$ is replaced by $\Phi_n$.*

Proof. We only treat the first inequality, the second is similar. Recall that by the LDP for the mass, we know that, for a fixed $n$,

$$\underline{\lim}_{\varepsilon\to0}\varepsilon\log\mathbb{P}(\mathbf{N}(u^{\varepsilon,0,n}(T))\ge R)\ge-\inf_{m>R}I^0_{N,n}(m), \tag{3.7}$$

where

$$I^0_{N,n}(m) = \frac{1}{2} \inf_{h \in L^2(0,T;L^2): \mathbf{N}(\mathbf{S}^{a,0,n}(h)(T))=m} \{\|h\|^2_{L^2(0,T;L^2)}\},$$

where $n$ means that $\Phi$ is replaced by $\Phi_n$ in the control equation.

We take $\tilde{R} > R$ and $h_{A_0}$ defined above. Though the stochastic equation is not defined when $\Phi = I$, the control map makes sense for any $h$ in $L^2(0,T;L^2)$. We denote it by $\mathbf{S}^{a,u_0}_{WN}$. By classical results (see [6], [10]), $\mathbf{S}^{a,0}_{WN}(h)$ is continuous with respect to $h \in C([0,T];L^2(\mathbb{R}))$ to $C([0,T];L^2(\mathbb{R}))$. Thanks to our assumptions, we deduce that

$$\mathbf{S}^{a,0,n}(h_{A_0}) = \mathbf{S}^{a,0}_{WN}(\Phi_n h_{A_0}) \to \mathbf{S}^{a,0}_{WN}(h_{A_0})$$

in $C([0,T];L^2(\mathbb{R}))$. In particular, since $N(\mathbf{S}^{a,0}_{WN}(h_{A_0})(T)) = \tilde{R}$, for $n$ large enough, we have $N(\mathbf{S}^{a,0,n}(h_{A_0})(T)) \geq R$ and the infimum in (3.7) is larger than

$$\frac{1}{2}\|h_{A_0}\|^2_{L^2(0,T;L^2)} = \frac{1}{18}(12+\pi^2)\frac{\tilde{R}}{T}.$$

We conclude since $\tilde{R}$ can be as close as we want to $R$.

Note that the result in Proposition 3.1 depends on $\Phi$ only through its norm as a bounded operator in $L^2$. It is not difficult to see that there exists sequences of operators $(\Phi_n)_{n\in\mathbb{N}}$ satisfying the assumptions of Proposition 3.2, that is, which are Hilbert–Schmidt from $L^2$ to $L^2$ and $\Phi_n$ approximates the identity on the good set of controls, and are uniformly bounded as operators on $L^2$ by a constant independent on $T$. For such sequences of operators, the upper and lower bounds given above agree up to constants in their behavior in large $T$. □

It is obtained in [19], for the ideal white noise and using the heuristic arguments recalled in the Introduction, that the probability density function of the amplitude of the pulse at coordinate $T$ when the initial datum is zero is asymptotically that of an exponential law of parameter $\varepsilon T/2$. The amplitude is a constant times the mass for the modulated soliton solutions considered [19]. Integrating this density over $[R/2,\infty)$ and taking into account the different normalization, we obtain $\lim_{\varepsilon\to 0} \varepsilon \log \mathbb{P}(\mathbf{N}(u^{\varepsilon,0}(T)) \geq R) = -R/T$. It is in between our two bounds and very close to our lower bound. A surprising fact is that we obtain our result by parameterizing only the amplitude, whereas in [19] a much more general parametrization is used. Both bounds exhibit the right behavior in $R$ and $T$. Moreover, the order in $R$ confirms physical and numerical results that the law is not Gaussian. On a log scale the order in $R$ is that of tails of an exponential law. In such a case the Gaussian approximation leads to incorrect tails and error estimates.

Let us now comment on our results in the case of a soliton profile as initial datum. In [19], the error probability when the size of the measurement window is of the order of the coordinate $T$ is obtained. It is given for some constant $c(R)$ by $\lim_{\varepsilon\to 0}\varepsilon\log\mathbb{P}(\mathbf{N}(u^{\varepsilon,\Psi}(T))-4\leq -R)=-c(R)/T$. It exhibits the same behavior in $T$ as in our calculations. The discussion on the behavior with respect to $R$ is less clear. Our bounds are not of the same order. In [14, 31] the PDF of the mass at coordinate $T$ for a soliton profile as initial datum is not Gaussian. The numerical simulations in [33] relying on the ansatz approximation also exhibit a significant difference between the log of the tails of the amplitude and that of a Gaussian law. Our lower bound indicates that again the tails are thicker than Gaussian tails. Thus, we give a rigorous proof of the fact that a Gaussian approximation is incorrect.

Finally, it is natural to obtain that the tails of the mass are increasing functions of $T$ since the higher is $T$, the less energy is needed to form a signal whose mass gets above a fixed threshold at $T$. Replacing above by under, the same holds in the case of a soliton as initial datum.

Remark 3.3. The $\mathrm{H}^1$ setting is not required here. We could as well work with $\mathrm{L}^2$ solutions and a LDP in $\mathrm{L}^2$. However, it is required to work in $\mathrm{H}^1$ for the study of the arrival time below.

We now estimate the tails of the arrival time. As for the mass, the rate is hard to handle since it involves an optimal control problem for controlled NLS equations. We again deduce the asymptotic of the tails from the LDP looking at upper and lower bounds. We consider that the initial datum is $\Psi_A^0$ since only in this case the timing jitter might be troublesome.

Let us begin with an upper bound. It is deduced from the equation of motion of the arrival time in the controlled NLS equation (2.10).

Proposition 3.4. *For every positive $T$, $A$ and $R$ and every operator $\Phi$ in $\mathcal{L}_2(\mathrm{L}^2,\Sigma)$, the following inequality holds:*

$$\begin{aligned}&\overline{\lim_{\varepsilon\to 0}}\,\varepsilon\log\mathbb{P}(\mathbf{Y}(u^{\varepsilon,\Psi_A^0}(T))\geq R)\\&\qquad\leq -\frac{R^2}{8T(2T+1)^2(4A+R/(2T+1))\|\Phi\|^2_{\mathcal{L}_c(\mathrm{L}^2,\Sigma)}}.\end{aligned}$$

Proof. Differentiating the right hand side of (2.10) with respect to time and replacing the time derivative of the solution with the corresponding terms of the equation we obtain

$$\frac{d}{dt}\mathbf{P}(\mathbf{S}^{a,\Psi_A^0}(h)(t))=4\mathfrak{Re}\int_{\mathbb{R}}\mathbf{S}^{a,\Psi_A^0}(h)(t,x)(\partial_x\overline{\Phi h})(t,x)\,dx.$$

Indeed, by successive integration by parts, all terms cancel besides the one involving the forcing term. Since $\mathbf{Y}(\Psi_A^0)=0$ and $\mathbf{P}(\Psi_A^0)=0$, thanks to (2.10), we obtain the identity

$$\mathbf{Y}(\mathbf{S}^{a,\Psi_A^0}(h)(t)) = 4\mathfrak{Re}\left(\int_0^t\int_0^s\int_{\mathbb{R}}\overline{\mathbf{S}^{a,\Psi_A^0}(h)(\sigma,x)}(\partial_x\Phi h)(\sigma,x)\,dx\,d\sigma\,ds\right)$$
$$-2\mathfrak{Re}\left(i\int_0^t\int_{\mathbb{R}}x\overline{\mathbf{S}^{a,\Psi_A^0}(h)(s,x)}(\Phi h)(s,x)\,dx\,ds\right).$$

From this identity it follows that the controls $h$ in the minimizing set of the LDP applied to the event we consider necessarily satisfy

$$R\le \mathbf{Y}(\mathbf{S}^{a,\Psi_A^0}(h)(T))\le 4T\|\Phi\|_{\mathcal{L}_c(\mathrm{L}^2,\mathrm{H}^1)}\|h\|_{\mathrm{L}^2(0,T;\mathrm{L}^2)}\|\mathbf{S}^{a,\Psi_A^0}(h)\|_{\mathrm{L}^2(0,T;\mathrm{L}^2)}$$
$$+2\|\Phi\|_{\mathcal{L}_c(\mathrm{L}^2,\Sigma)}\|h\|_{\mathrm{L}^2(0,T;\mathrm{L}^2)}\|\mathbf{S}^{a,\Psi_A^0}(h)\|_{\mathrm{L}^2(0,T;\mathrm{L}^2)}.$$

Moreover, arguing as in the proof of Proposition 3.1,

$$\|\mathbf{S}^{a,\Psi_A^0}(h)\|_{\mathrm{L}^2(0,T;\mathrm{L}^2)}\le T\|\Phi\|_{\mathcal{L}_c(\mathrm{L}^2,\mathrm{L}^2)}\|h\|_{\mathrm{L}^2(0,T;\mathrm{L}^2)}$$
$$\times\left(1+\sqrt{1+4A/(T\|\Phi\|^2_{\mathcal{L}_c(\mathrm{L}^2,\mathrm{L}^2)}\|h\|^2_{\mathrm{L}^2(0,T;\mathrm{L}^2)})}\right).$$

A lower bound on $(1/2)\|h\|^2_{\mathrm{L}^2(0,T;\mathrm{L}^2)}$ follows easily since $x\mapsto x(1+\sqrt{1+4/x})$ is increasing on $\mathbb{R}_+^*$. The result follows. □

A lower bound is obtained considering controls suggested at the end of Section 2 and minimizing on the smaller set of controls. We define the following set of control for $A,T$ positive and $D$ a subset of $(0,\infty)$:

$$\mathcal{H}^D_{A,T}=\Big\{h\in\mathrm{L}^2(0,T;\mathrm{L}^2),h(t,x)=\lambda(t)\left(x-2\int_0^t\int_0^s\lambda(\tau)\,d\tau\,ds\right)\tilde{\Psi}_{A,\lambda}(t,x),$$
$$\text{with }\lambda(t)=3\tilde{R}(T-t)/(8AT^3),\tilde{R}\in D\Big\}.$$

Proposition 3.5. *Let $T$, $A$ and $R$ be positive. Assume that, for a dense subset $D$ of $[R,R+1]$, $(\Phi_n)_{n\in\mathbb{N}}$ is a sequence of operators in $\mathcal{L}_2(\mathrm{L}^2,\Sigma)$ such that for any $h$ in $\mathcal{H}^D_{T,A}$, $\Phi_n h$ converges to $h$ in $\mathrm{L}^1(0,T;\Sigma)$. Then we have the following inequality where the exponent $n$ is there to recall that $\Phi$ is replaced by $\Phi_n$:*

$$\varliminf_{n\to\infty,\varepsilon\to0}\varepsilon\log\mathbb{P}(\mathbf{Y}(u^{\varepsilon,\Psi_A^0,n}(T))\ge R)\ge-\pi^2R^2/(128T^3A^3).$$

Proof. We proceed as for Proposition 3.2. By the LDP for the arrival time $\mathbf{Y}$, we know that for a fixed $n$ a lower bound is given by

$$-\inf_{y>R}I^{\Psi_A^0}_{Y,n}(y),$$

where

$$I^{\Psi^0_A}_{Y,n}(y) = \frac{1}{2} \inf_{h\in \mathrm{L}^2(0,T;\mathrm{L}^2):\mathbf{Y}(\mathbf{S}^{a,\Psi^0_A,n}(h)(T))=y} \{\|h\|^2_{\mathrm{L}^2(0,T;\mathrm{L}^2)}\}.$$

Again, $n$ is there to recall that in the control equation, $\Phi$ is replaced by $\Phi_n$. To minorize this quantity, we again first treat the case $\Phi = I$ and denote by $\mathbf{S}^{a,\Psi^0_A}_{WN}$ the control map when $\Phi = I$.

It is not difficult to see that $\mathbf{S}^{a,\Psi^0_A}_{WN}(h)$ belongs to $\mathrm{L}^\infty([0,T];\Sigma)$ when $h$ belong to $\mathrm{L}^1(0,T;\Sigma)$. Moreover, the norm of $\mathbf{S}^{a,\Psi^0_A}_{WN}(h)$ in $\mathrm{L}^\infty([0,T];\Sigma)$ is bounded in terms of the norm of $h$ in $\mathrm{L}^1(0,T;\Sigma)$. A standard argument to prove this is to compute the second derivative with respect to time of the variance $\mathbf{V}(u) = \int_{\mathbb{R}} x^2|u(t,x)|^2\,dx$ with $u = \mathbf{S}^{a,\Psi^0_A}_{WN}(h)$. Using the argument detailed in the proof of the LDP in Section 5 below, this implies that, for each $t$, the mapping $h \to \mathbf{S}^{a,\Psi^0_A}_{WN}(h)(t)$ is continuous from $\mathrm{L}^1(0,T;\Sigma)$ to $\Sigma^{1/2}$. Therefore, for $h \in \mathcal{H}^D_{T,A}$,

$$\begin{aligned}\mathbf{Y}(\mathbf{S}^{a,\Psi^0_A,n}(h)(T)) &= \mathbf{Y}(\mathbf{S}^{a,\Psi^0_A}_{WN}(\Phi_n h)(T)) \\ &\to \mathbf{Y}(\mathbf{S}^{a,\Psi^0_A}_{WN}(h)(T)) \qquad \text{when } n\to\infty.\end{aligned} \tag{3.8}$$

Proceeding as above, we are thus led to find a control $h$ with minimum energy verifying the constraint $\mathbf{Y}(\mathbf{S}^{a,\Psi^0_A}_{WN}(h)(T)) \ge \tilde{R}$ for some $\tilde{R} > R$.

We search this control in the set $\tilde{\mathcal{H}}_{T,A}$ defined as $\mathcal{H}_{T,A}$, but where $\lambda$ is not specified and only assumed to belong to $\mathrm{L}^2(0,T;\mathbb{R})$:

$$\tilde{\mathcal{H}}_{T,A} = \Big\{ h \in \mathrm{L}^2(0,T;\mathrm{L}^2), \\ h(t,x) = \lambda(t)\Big(x - 2\int_0^t\int_0^s \lambda(\tau)\,d\tau\,ds\Big)\tilde{\Psi}_{A,\lambda}(t,x), \lambda\in \mathrm{L}^2(0,T;\mathbb{R})\Big\}.$$

We have seen at the end of Section 2 that $\mathbf{S}^{a,\Psi^0_A}_{WN}(h) = \tilde{\Psi}_{A,\lambda}$ for $h \in \tilde{\mathcal{H}}_{T,A}$. Also, an easy computation gives $\mathbf{Y}(\tilde{\Psi}_{A,\lambda}(T)) = 8A\int_0^T\int_0^t \lambda(s)\,ds\,dt$ and, for $h \in \tilde{\mathcal{H}}_{T,A}$, $\|h\|^2_{\mathrm{L}^2(0,T;\mathrm{L}^2)} = (\pi^2/(3A))\int_0^T \lambda^2(s)\,ds$. We deduce

$$\begin{aligned}&\inf_{h\in \mathrm{L}^2(0,T;\mathrm{L}^2):\mathbf{Y}(\mathbf{S}^{a,\Psi^0_A}_{WN}(h)(T))\ge\tilde{R}} \|h\|^2_{\mathrm{L}^2(0,T;\mathrm{L}^2)} \\ &\quad\le \inf_{h\in\tilde{\mathcal{H}}_{T,A}:\mathbf{Y}(\mathbf{S}^{a,\Psi^0_A}_{WN}(h)(T))\ge\tilde{R}} \|h\|^2_{\mathrm{L}^2(0,T;\mathrm{L}^2)} \\ &\quad= \inf_{\lambda\in \mathrm{L}^2(0,T;\mathbb{R}):\int_0^T\int_0^t\lambda(s)\,ds\,dt\ge\tilde{R}/(8A)} (\pi^2/(3A))\int_0^T \lambda^2(t)\,dt.\end{aligned}$$

Note that the constraint $\int_0^T\int_0^t \lambda(s)\,ds\,dt \geq \tilde{R}/8A$ is not a usual boundary condition in the calculus of variations. We therefore simply try to find a good guess that leads to a lower bound of the the same order in the parameters as the upper bound. We use the quantity $\mathcal{L}_{T,A,\tilde{R}}(\lambda)$ defined by

$$\mathcal{L}_{T,A,\tilde{R}}(\lambda) = (\pi^2/(3A))\int_0^T \lambda^2(t)\,dt - \gamma\int_0^T\int_0^t \lambda(s)\,ds\,dt,$$

where $\gamma$ belongs to $\mathbb{R}$. We then impose that our guess $\lambda^*_{T,A,\tilde{R}}$ is a critical point of $\mathcal{L}_{T,A,\tilde{R}}(\lambda)$ and that it satisfies the constraint $\int_0^T\int_0^t \lambda(s)\,ds\,dt = \tilde{R}/(8A)$. We obtain

$$\lambda^*_{T,A,\tilde{R}}(t) = 3\tilde{R}(T-t)/(8AT^3).$$

We do not claim that the minimization problem is solved, we simply write

$$\begin{aligned}\inf_{\lambda\in \mathrm{L}^2(0,T;\mathbb{R}):\int_0^T\int_0^t \lambda(s)\,ds\,dt\geq \tilde{R}/(8A)} &(\pi^2/(3A))\int_0^T \lambda^2(t)\,dt \\ &\leq (\pi^2/(3A))\int_0^T (\lambda^*_{T,A,\tilde{R}}(t))^2\,dt = \pi^2\tilde{R}^2/(64A^3T^3).\end{aligned}$$

Let us set

$$h^*_{\tilde{R}}(t,x) = \lambda^*_{T,A,\tilde{R}}(t)\left(x - 2\int_0^t\int_0^s \lambda^*_{T,A,\tilde{R}}(\tau)\,d\tau\,ds\right)\tilde{\Psi}_{A,\lambda^*_{T,A,\tilde{R}}}(t,x).$$

By (3.8), we have, for $\tilde{R}\in D$,

$$\mathbf{Y}(\mathbf{S}^{a,\Psi^0_A,n}(h^*_{\tilde{R}})(T)) \to \mathbf{Y}(\mathbf{S}^{a,\Psi^0_A}_{WN}(h^*_{\tilde{R}})(T)) \qquad \text{when } n\to\infty.$$

Therefore, for $n$ large enough,

$$\mathbf{Y}(\mathbf{S}^{a,\Psi^0_A,n}(h^*_{\tilde{R}})(T)) > R.$$

We deduce

$$\inf_{x>R} I^{\Psi^0_A}_{Y,n}(x) \leq \pi^2\tilde{R}^2/(64A^3T^3)$$

and take the $\underline{\lim}$ in $n$ in the lower bound. Since this is true for $\tilde{R}$ in a dense subset of $[R, R+1]$, we deduce the result. □

The upper and lower bounds given in Propositions 3.4 and 3.5 are in perfect agreement in their behavior with respect to $R$ and to $T$ when $T$ is large. Indeed, as $T$ is large compared to $R$, the upper bound in Proposition 3.4 is of the order of $-R^2/(128T^3A\|\Phi\|^2_{\mathcal{L}_c(\mathrm{L}^2,\Sigma)})$. However, we have to be careful before doing such a comparison. Indeed, unlike for the mass, it does not seem that there exists a sequence of operators $(\Phi_n)_{n\in\mathbb{N}}$ satisfying the assumptions of Proposition 3.5 and such that $\|\Phi_n\|_{\mathcal{L}_c(\mathrm{L}^2,\Sigma)}$ is bounded uniformly in $n$. This explains why the behavior in $A$ in the lower and upper bound seem contradictory for large $A$. We believe, however, that there exists a sequence satisfying the assumptions of Proposition 3.5 and such that $\|\Phi_n\|_{\mathcal{L}_c(\mathrm{L}^2,\Sigma)}$ is bounded independently with respect to $R$ and $T$, but not to $A$.

It is, however, possible to obtain bounds that match with respect to their order in $A$. Unlike the framework of Proposition 3.5, suppose we consider the sequence of operators $\Phi_n = (I - \Delta + |x|^2 I + \frac{1}{n}(-\Delta + |x|^2 I)^k)^{-1/2}$ such that $\Phi_n h \to \Phi h$ for any $h$ and $\Phi = (I - \Delta + |x|^2 I)^{-1/2}$. We may prove that for sufficiently large $k$ $\Phi_n$ are $\Sigma$-valued Hilbert–Schmidt operators. Also then $\|\Phi_n\|_{\mathcal{L}_c(\mathrm{L}^2,\Sigma)} \leq \|\Phi\|_{\mathcal{L}_c(\mathrm{L}^2,\Sigma)} = 1$ and, thus, $\Phi_n$ are bounded uniformly in $n$. We argue as in the above using that

$$\mathbf{Y}(\mathbf{S}_{WN}^{a,\Psi_A^0}(\Phi_n h_{\tilde{R}}^*)(T)) \to \mathbf{Y}(\mathbf{S}_{WN}^{a,\Psi_A^0}(\Phi h_{\tilde{R}}^*)(T)) \qquad \text{when } n \to \infty.$$

Then for $n$ large enough, the lower bound is given by $-\|h_{\tilde{R}}^*\|_{\mathrm{L}^2(0,T;\Sigma)}^2/2$. Then for large $A$, that is, very localized initial pulse allowing theoretically higher transmission rate, the order up to a multiplicative constant is now that of the square of the norm of the gradient. It is thus now $-\tilde{R}/(T^3 A)$ which matches the order of the upper bound.

Let us now compare our result with the results obtained in the physics literature. First, we note that we obtain that on a log scale the tails are equivalent to Gaussian tails. This is indeed the kind of result obtained by arguments from the physical theory of perturbation of solitons. We are missing the pre exponential factors to conclude whether or not the tails are Gaussian. Sharp large deviations could allow to obtain these factors.

Now, suppose the law were indeed Gaussian, then the asymptotic of the tails may be written in terms of the variance. By doing so, we find that the variance of the timing jitter is of the order $T^3$. It agrees perfectly with the initial results of [24]. Also, the order in $A$—for the lower bound—and $T$ seems to agree perfectly with the orders of the contribution of the additive noise to the variance of the timing jitter in equation (3.18) in [17]. However it is not clear that we wish to obtain the $1/A^3$ order since in [24, 28], where the model is instead a juxtaposition of deterministic evolutions with random initial data in between amplifiers, the order in $A$ is $-c/A$. It is what we obtain above when we assume for consistency of the assumptions that in the limit the noise remains localized.

We end this section noticing that our result confirms the fact that, in the presence of additive noise, the timing jitter is more troublesome than the fluctuation of the mass when we consider the problem of losing a signal. Indeed, for $A = 1$, we have found that the tails of the arrival time are the order of $\exp(-c_1(R)/(\varepsilon T^3))$, while that of the mass are of the order of $\exp(-c_2(R)/(\varepsilon T))$, which is clearly negligible compared to the first for large $T$. In other words the tails of the arrival time are much thicker than that of the mass, implying much more frequent large fluctuations of the arrival time than of the mass. Error in soliton transmission is much more likely to be due to timing jitter rather than decay of the mass of the pulse. Recall that $T$ is the length of a fiber optical line and is thus assumed to be very large. This result is called the Gordon–Haus effect in the physics literature.

REMARK 3.6. From an engineering point of view, it is possible to exponentially reduce the probability of undesired deviations of the arrival time by introducing inline control elements; see, for example, [19]. We could also use ideas given in [36] and optimize on such external fields for a limited cost or penalty functional. The new optimal control problem requires then double optimization.

REMARK 3.7. Note that the methodology developed herein could probably be applied to the determination of the small noise asymptotic of the tails of the position of an isolated vortex, defined by $\oint \nabla \arg u(t,x) \cdot d\mathbf{l}$, in Bose condensates or superfluid Helium as in [34]. There the physical perturbation approach along with the Fokker–Planck equation are used. The small noise acts as the small temperature.

**4. Tails of the arrival time in the multiplicative case.** In the case of the multiplicative noise, the mass is a conserved quantity and we restrict our attention to the study of the law of the arrival time of the pulse when the initial datum is the soliton profile $\Psi_A^0$.

Again, let us begin with upper bounds. They are obtained from an equation for the motion of the arrival time in the controlled NLS equation.

From relation (2.9) and integration by parts, we obtain the equation in [38],

$$\frac{d^2}{dt^2}\mathbf{Y}(\mathbf{S}^{m,\Psi_A^0}(h)(t)) = 2\int_{\mathbb{R}} |\mathbf{S}^{m,\Psi_A^0}(h)(t,x)|^2 (\partial_x \Phi h)(t,x)\,dx. \tag{4.1}$$

We may thus deduce the next proposition.

PROPOSITION 4.1. *For every positive $T$, $A$ and $R$ and every operator $\Phi$ in $\mathcal{L}_2(\mathrm{L}^2,\mathrm{H}^s(\mathbb{R},\mathbb{R}))$, where $s>3/2$ the following inequality holds:*

$$\overline{\lim_{\varepsilon\to 0}}\, \varepsilon \log \mathbb{P}(\mathbf{Y}(u^{\varepsilon,\Psi_A^0}(T)) \geq R) \leq -\left(\frac{3}{16}\right)^2 \frac{R^2}{2A^2T^3\|\Phi\|^2_{\mathcal{L}_c(\mathrm{L}^2,\mathrm{W}^{1,\infty}(\mathbb{R},\mathbb{R}))}}.$$

PROOF. From equation (4.1), the fact that $\frac{d}{dt}\mathbf{Y}(\mathbf{S}^{m,\Psi_A^0}(h))|_{t=0} = \mathbf{P}(\Psi_A^0) = 0$, that for such values of $s$ the Sobolev injection of $\mathrm{H}^s(\mathbb{R},\mathbb{R})$ into $\mathrm{W}^{1,\infty}(\mathbb{R},\mathbb{R})$

is continuous (see [18]), and that the mass is conserved and thus remains equal to 4, we obtain

$$\frac{d}{dt}\mathbf{Y}(\mathbf{S}^{m,\Psi_A^0}(h)(t)) \le 8A\|\Phi\|_{\mathcal{L}_c(\mathrm{L}^2,\mathrm{W}^{1,\infty}(\mathbb{R},\mathbb{R}))}\|h\|_{\mathrm{L}^1(0,t;\mathrm{L}^2)}$$
$$\le 8A\sqrt{t}\|\Phi\|_{\mathcal{L}_c(\mathrm{L}^2,\mathrm{W}^{1,\infty}(\mathbb{R},\mathbb{R}))}\|h\|_{\mathrm{L}^2(0,T;\mathrm{L}^2)}.$$

Then, since $\mathbf{Y}(\Psi_A^0)=0$, we obtain integrating the above inequality that

$$R \le \mathbf{Y}(\mathbf{S}^{m,\Psi_A^0}(h)(T)) \le (16AT^{3/2}/3)\|\Phi\|_{\mathcal{L}_c(\mathrm{L}^2,\mathrm{W}^{1,\infty}(\mathbb{R},\mathbb{R}))}\|h\|_{\mathrm{L}^2(0,T;\mathrm{L}^2)}$$

and the conclusion follows. □

Let us now consider lower bounds. We need to find controls which have the desired effect on the arrival time. We have seen that, in the additive case, good controls are given by functions in $\mathcal{H}_{A,T}^D$. Recalling the transformations on the equation made at the end of Section 2, we can equivalently take controls of the form $\lambda(t)x\Psi_{A,\lambda}$ which correspond to the solution $\Psi_{A,\lambda}$. Thus, in the multiplicative case, a good control is given by $h(t,x)=\lambda(t)x$. Unfortunately these controls do not belong to the range of $\Phi$ nor to $\mathrm{L}^2(0,T;\mathrm{L}^2)$ and are not admissible. We have not been able to justify the choice of such controls by an approximation argument.

We therefore impose a new assumption that $\Phi$ takes its values in $\mathrm{H}^s(\mathbb{R},\mathbb{R})\oplus x\mathbb{R}$. In other words, we consider the slightly different equation

$$\begin{aligned} id\tilde{u}^{\varepsilon,u_0} &= (\Delta\tilde{u}^{\varepsilon,u_0} + |\tilde{u}^{\varepsilon,u_0}|^2\tilde{u}^{\varepsilon,u_0})\,dt \\ &\quad + \tilde{u}^{\varepsilon,u_0}\circ\sqrt{\varepsilon}\,dW(t) + \sqrt{\varepsilon}x\tilde{u}^{\varepsilon,u_0}\circ d\beta(t), \end{aligned} \tag{4.2}$$

where $\beta$ is a standard Brownian motion independent of $W$ and the corresponding controlled PDE

$$i\frac{d}{dt}\tilde{S}^{u_0}(h_1,h_2) = \Delta\tilde{S}^{u_0}(h_1,h_2) + |\tilde{S}^{u_0}(h_1,h_2)|^2\tilde{S}^{u_0}(h_1,h_2)$$
$$+ \tilde{S}^{u_0}(h_1,h_2)\Phi h_1 + x\tilde{S}^{u_0}(h_1,h_2)h_2,$$

where $h_1$ belongs to $\mathrm{L}^2(0,T;\mathrm{L}^2)$ and $h_2$ belongs to $\mathrm{L}^2(0,T;\mathbb{R})$, the initial datum is $u_0$ and in the sequel $u_0=\Psi_A^0$. We may guess by successive applications of the Itô formula, multiplying $\tilde{u}^{\varepsilon,u_0}$ by the random phase term $\exp(ix\sqrt{\varepsilon}\beta(t))$, and similar transformations as in Section 2 (stochastic gauge transform, stochastic methods of characteristics, ...) that we should consider the function

$$\exp\Big(ix\sqrt{\varepsilon}\beta(t) - i\varepsilon\int_0^t \beta^2(s)\,ds\Big)\tilde{u}^{\varepsilon,u_0}\Big(t, x+2\sqrt{\varepsilon}\int_0^t \beta(s)\,ds\Big).$$

It indeed satisfies equation (2.3) with the same initial datum. We deduce that

$$\tilde{u}^{\varepsilon,u_0}(t,x) = \exp\Big(-ix\sqrt{\varepsilon}\beta(t) + i\varepsilon\int_0^t \beta^2(s)\,ds + 2i\varepsilon\beta(t)\int_0^t \beta(s)\,ds\Big) \times u^{\varepsilon,u_0}\Big(t, x - 2\sqrt{\varepsilon}\int_0^t \beta(s)\,ds\Big).$$

A similar computation shows that

$$\tilde{S}^{u_0}(h_1,h_2)(t,x) = \exp\Big(-ix\sqrt{\varepsilon}\int_0^t h_2(s)\,ds + i\int_0^t\Big(\int_0^s h_2(u)\,du\Big)^2 ds + 2i\int_0^t h_2(s)\,ds\int_0^t\int_0^s h_2(u)\,du\,ds\Big) \times \mathbf{S}^{m,u_0}(h_1)\Big(t, x - 2\int_0^t\int_0^s h_2(u)\,du\Big).$$

The functions $\tilde{u}^{\varepsilon,u_0}$ and $\tilde{S}^{u_0}(h_1,h_2)$ are well-defined functions of $\mathrm{L}^2(0,T;\Sigma)$ and we may compute their arrival times. We obtain a lower bound of the asymptotic of the tails of the arrival time of the new solutions.

Proposition 4.2. *For every positive $T$, $A$ and $R$ and every operator $\Phi$ in $\mathcal{L}_2(\mathrm{L}^2,\mathrm{H}^s(\mathbb{R},\mathbb{R}))$ where $s>3/2$, the following inequality holds:*

$$\varliminf_{\varepsilon\to 0}\varepsilon\log\mathbb{P}(\mathbf{Y}(\tilde{u}^{\varepsilon,\Psi_A^0}(T)) \ge R) \ge -3R^2/(128A^2T^3).$$

Proof. Consider the mapping $F$ from $\mathrm{C}([0,T];\Sigma^{1/2})\times\mathrm{C}([0,T];\mathbb{R})$ into $\mathbb{R}$ such that

$$F(u,b) = \int_{\mathbb{R}} |x|\Big|u\Big(T, x - 2\int_0^T b(s)\,ds\Big)\Big|^2 dx.$$

Take $u$ and $u'$ in $\mathrm{C}([0,T];\Sigma^{1/2})$ and $b$ and $b'$ in $\mathrm{C}([0,T];\mathbb{R})$, then by the triangle and inverse triangle inequalities and the change of variables, we obtain

$$\begin{aligned}
&|F(u,b) - F(u',b')| \\
&\quad\le \int_{\mathbb{R}}\bigg|\Big|x + 2\int_0^T b(s)\,ds\Big| - \Big|x + 2\int_0^T b'(s)\,ds\Big|\bigg|\,|u(T,x)|^2\,dx \\
&\qquad + \bigg|\int_{\mathbb{R}}\Big|x + 2\int_0^T b'(s)\,ds\Big|(|u(T,x)|^2 - |u'(T,x)|^2)\,dx\bigg| \\
&\quad\le 2\Big|\int_0^T b(s)\,ds - \int_0^T b'(s)\,ds\Big|\int_{\mathbb{R}}|u(T,x)|^2\,dx
\end{aligned}$$

$$+\int_{\mathbb{R}}|x|\,||u(T,x)|-|u'(T,x)||(|u(T,x)|+|u'(T,x)|)\,dx$$
$$+2\left|\int_0^T b'(s)\,ds\right|\int_{\mathbb{R}}||u(T,x)|-|u'(T,x)||(|u(T,x)|+|u'(T,x)|)\,dx.$$

We conclude from the inverse triangle and Hölder inequalities that $F$ is continuous. We may then push forward by the contraction principle the LDP for the paths of $u^{\varepsilon,\Psi_A^0}$ and of $\sqrt{\varepsilon}\beta$ by the mapping $F$ using a slight modification of the result of exercise 4.2.7 of [9] and obtain a LDP for the laws of $\mathbf{Y}(\tilde{u}^{\varepsilon,\Psi_A^0}(T))$ which is that of $F(u^{\varepsilon,\Psi_A^0},\sqrt{\varepsilon}\beta)$ of speed $\varepsilon$ and good rate function defined as a function of the rate function of the original solutions and of the rate function $I_\beta$ of the sample path LDP for $\sqrt{\varepsilon}\beta$:

$$\begin{aligned}\tilde{I}_Y^{\Psi_A^0}(x)&=\inf_{(u,b):F(u,b)=x}(I^{u_0}(u)+I_\beta(b))\\&\le\tfrac12\inf_{(h_1,h_2):F(\mathbf{S}^{m,\Psi_A^0}(h_1),\int_0^\cdot h_2(s)ds)=x}\{\|h_1\|^2_{\mathrm{L}^2(0,T;\mathrm{L}^2)}+\|h_2\|^2_{\mathrm{L}^2(0,T;\mathbb{R})}\}\\&\le\tfrac12\inf_{(h_1,h_2):\mathbf{Y}(\tilde{\mathbf{S}}^{\Psi_A^0}(h_1,h_2)(T))=x}\{\|h_1\|^2_{\mathrm{L}^2(0,T;\mathrm{L}^2)}+\|h_2\|^2_{\mathrm{L}^2(0,T;\mathbb{R})}\}.\end{aligned}$$

Thus, considering solely controls of the from $(0,h_2)$, we minimize in $h_2$ for $\gamma$ in $\mathbb{R}$,

$$\int_0^T h_2^2(t)\,dt-\gamma\int_0^T\int_0^t h_2(s)\,ds,$$

where we impose that

$$\mathbf{Y}(\Psi_{A,h_2}(T))=8A\int_0^T\int_0^t h_2(s)\,ds=\tilde{R}>R.$$

The conclusion follows. □

Remark 4.3. We may check that $\mathbf{Y}(u^{\varepsilon,\Psi_A^0})=\mathbf{Y}(\tilde{u}^{\varepsilon,\Psi_A^0})-8\sqrt{\varepsilon}\int_0^T\beta(s)\,ds$ and that $\int_0^T\beta(s)\,ds$ is a centered Gaussian random variable with variance $T^3/3$.

The corresponding upper bound for this modified stochastic NLS equation is

$$\overline{\lim_{\varepsilon\to0}}\,\varepsilon\log\mathbb{P}(\mathbf{Y}(\tilde{u}^{\varepsilon,\Psi_A^0}(T))\ge R)\le-(3/16)^2\frac{R^2}{A^2T^3(\|\Phi\|^2_{\mathcal{L}_c(\mathrm{L}^2,\mathrm{W}^{1,\infty}(\mathbb{R},\mathbb{R}))}\vee1)}.$$

Note that the lower bound does not require to consider a sequence of operators $(\Phi_n)_{n\in\mathbb{N}}$ and we may indeed compare the upper and lower bounds.

They are of the same order in $T$ and in $A$. Note also that, as in the additive case, we obtain that on a log scale the tails are equivalently that of Gaussian tails. Also, our tails are of the order in $T$ that we expect from the contribution of the multiplicative noise to the variance of the timing jitter in equation (3.18) in [17].

However, concerning the amplitude, it is not of the order of $-c/A^4$ as we would expect from [17]. This is probably due to the fact that we have considered a colored noise with a term $x\frac{d}{dt}\beta$ that grows linearly in time (the $x$ variable).

**5. Proof of Theorem 2.1.** We herein denote the variance of $f$ in $\Sigma$ by $\mathbf{V}(f)=\int_{\mathbb{R}}|x|^2|f(x)|^2\,dx$.

Let us start with the additive case. We denote by $v^{u_0}(z)$ the solution of

$$\begin{cases} i\dfrac{dv}{dt}=\Delta v+\lambda|v-iz|^{2\sigma}(v-iz),\\ u(0)=u_0\in\Sigma,\end{cases}$$

where $z$ belongs to $X^{(T,2\sigma+2,\Sigma)}=\mathrm{C}([0,T];\Sigma)\cap\mathrm{L}^r(0,T;\mathrm{W}^{1,2\sigma+2})$ and $r$ is such that $2/r=1/2-1/(2\sigma+2)$. We also denote by $\mathcal{G}^{u_0}$ the mapping

$$z\mapsto v^{u_0}(z)-iz.$$

Note that $u^{\varepsilon,u_0}=\mathcal{G}^{u_0}(\sqrt{\varepsilon}Z)$, where $Z$ is the stochastic convolution defined by $Z(t)=\int_0^t U(t-s)\,dW(s)$.

We can check from similar arguments (as those of the proof of Proposition 1 in [22]) that the stochastic convolution is a $X^{(T,2\sigma+2,\Sigma)}$ random variable whose law $\mu^Z$ is a centered Gaussian measure. Let $z$ belong to $X^{(T,2\sigma+2,\Sigma)}$, take $s<t<T$, the triangle along with the Hölder inequalities, then compute

$$\begin{aligned}&\left|\int_{\mathbb{R}}|x|(|\mathcal{G}^{u_0}(z)(t,x)|^2-|\mathcal{G}^{u_0}(z)(s,x)|^2)\,dx\right|\\
&\qquad\le\int_{\mathbb{R}}|x|(|\mathcal{G}^{u_0}(z)(t,x)|+|\mathcal{G}^{u_0}(z)(s,x)|)|(|\mathcal{G}^{u_0}(z)(t,x)|-|\mathcal{G}^{u_0}(z)(s,x)|)|\,dx\\
&\qquad\le\|\mathcal{G}^{u_0}(z)(t)-\mathcal{G}^{u_0}(z)(s)\|_{\mathrm{L}^2}\sqrt{\mathbf{V}(|\mathcal{G}^{u_0}(z)(t)|+|\mathcal{G}^{u_0}(z)(s)|)}\\
&\qquad\le 2\sqrt{2}\|\mathcal{G}^{u_0}(z)(t)-\mathcal{G}^{u_0}(z)(s)\|_{\mathrm{L}^2}\\
&\qquad\quad\times\left(\sqrt{\mathbf{V}(v^{u_0}(z)(t))}+\sqrt{\mathbf{V}(v^{u_0}(z)(s))}+\sqrt{\mathbf{V}(z(t))}+\sqrt{\mathbf{V}(z(s))}\right).\end{aligned}$$

The application of the Gronwall inequality in the proof of Proposition 3.5 in [11], along with the Sobolev injection allow to prove that $\mathcal{G}^{u_0}(z)$ belongs to $\mathrm{C}([0,T];\Sigma^{1/2})$. The computation above also shows that $\mathcal{G}^{u_0}$ is continuous from $X^{(T,2\sigma+2,\Sigma)}$ to $\mathrm{C}([0,T];\Sigma^{1/2})$. The general result on LDP for Gaussian measures gives the LDP for the measures $\mu^{Z_\varepsilon}$, the direct images of $\mu^Z$ under the transformation $x \mapsto \sqrt{\varepsilon} x$ on $X^{(T,2\sigma+2,\Sigma)}$. We conclude with the contraction principle.

In the multiplicative case, it is also required to revisit the proof of the LDP in [23]. As mentioned in Section 2, we only emphasize the adaptations of the proof of [23] required to state a LDP in $\mathrm{C}([0,T];\Sigma^{1/2})$. Note that in the following $\Phi h$ is replaced by $\frac{\partial f}{\partial t}$, where $f$ belongs to $\mathrm{H}_0^1(0,T;\mathrm{H}^s(\mathbb{R},\mathbb{R}))$ which is the subspace of $\mathrm{C}([0,T];\mathrm{H}^s(\mathbb{R},\mathbb{R}))$ of functions zero at time 0, square integrable in time and with square integrable in time derivative. The control map is then denoted by $\tilde{\mathbf{S}}^{m,u_0}(f)$.

We may check using the above calculation and the fact that, for every $t \in [0,T]$, $\tilde{\mathbf{S}}^{m,u_0}(f)(t)$ belongs to $\Sigma$ that

$$\mathbf{V}(\tilde{\mathbf{S}}^{m,u_0}(f)(t)) \le (4\|\tilde{\mathbf{S}}^{m,u_0}(f)(t)\|^2_{\mathrm{C}([0,T];\mathrm{H}^1)} + \mathbf{V}(u_0))e^T;$$

(see the arguments of the proof of Proposition 3.2 in [12]) used for the control map, that the control map is continuous from the sets of levels of the rate function of the Wiener process less or equal to a positive constant, with the topology induced by that of $\mathrm{C}([0,T];\mathrm{H}^s(\mathbb{R},\mathbb{R}))$, to $\mathrm{C}([0,T];\Sigma^{1/2})$. The only difference in the proof of Proposition 4.1 in [23], the Azencott lemma (also called the Freidlin–Wentzell inequality or almost continuity of the Itô map) is in step 2. It is the reduction to estimates on the stochastic convolution. We use

$$\mathbf{V}(v^{\varepsilon,\tilde{u}_0}(t)) \le (4\|v^{\varepsilon,\tilde{u}_0}(t)\|^2_{\mathrm{C}([0,T];\mathrm{H}^1)} + \mathbf{V}(\tilde{u}_0))e^T$$

(see the proof of Proposition 3.2 in [12]), where $v^{\varepsilon,\tilde{u}_0}$ satisfies $v^{\varepsilon,\tilde{u}_0}(0) = \tilde{u}_0$ and

$$\begin{aligned} i dv^{\varepsilon,\tilde{u}_0} &= \left( \Delta v^{\varepsilon,\tilde{u}_0} + \lambda |v^{\varepsilon,\tilde{u}_0}|^{2\sigma} v^{\varepsilon,\tilde{u}_0} + \frac{\partial f}{\partial t} v^{\varepsilon,\tilde{u}_0} - (i\varepsilon/2) F_\Phi v^{\varepsilon,\tilde{u}_0} \right) dt \\ &\quad + \sqrt{\varepsilon} v^{\varepsilon,\tilde{u}_0}\, dW_\varepsilon, \end{aligned}$$

$f(\cdot) = \int_0^\cdot \Phi h(s)\, ds$, $W_\varepsilon(t) = W(t) - (1/\sqrt{\varepsilon}) \int_0^t \frac{\partial f}{\partial s}\, ds = W(t) - (1/\sqrt{\varepsilon}) \times \int_0^t \Phi h(s)\, ds$, $F_\Phi(x) = \sum_{j=1}^\infty (\Phi e_j(x))^2$ and $(e_j)_{j=1}^\infty$ is any complete orthonormal system of $\mathrm{L}^2$. The bound remains the same as in [12] because of the cancelation of the extra term in the application of the Itô formula and the cancelation of the Itô–Stratonovich correction with the second order Itô correction term when the Itô formula is applied to the truncated variance $V_r(v) = \int_{\mathbb{R}} \exp(-r|x|^2)|x|^2|v(x)|^2\, dx$.

**Acknowledgments.** Eric Gautier is also affiliated to CREST, Laboratoire de Statistique, 3 avenue Pierre Larousse, 92240 Malakoff, France.

IRMAR, ENS Cachan Bretagne, CNRS, UEB
av Robert Schuman
F-35170 Bruz
France
E-mail: arnaud.debussche@bretagne.ens-cachan.fr

IRMAR, ENS Cachan Bretagne, CNRS, CREST
av Robert Schuman
F-35170 Bruz
France
E-mail: gautier@ensae.fr